# Counting integral points of affine hypersurfaces

*Per Salberger*

ABSTRACT  We give uniform upper bounds for the number of integral points of bounded height on affine hypersurfaces, which generalise earlier results of Browning, Heath-Brown and the author.

## 1 Introduction

We shall in this note generalise results in [BHS], [$S_2$] and [$S_3$] and prove the following theorem about the density of integral points on affine hypersurfaces.

**Theorem 1**  *Let $f(x_1,...,x_n) \in \mathbf{Q}[x_1,..., x_n]$ be an irreducible polynomial of degree d in $n \geq 3$ variables such that f cannot be expressed as a polynomial in two linear forms in $(x_1,..., x_n)$ and let $n(f;B)$ be the number of integral n-tuples $(a_1,...,a_n)$ in $[-B, B]^n$ with $f(a_1,..., a_n)=0$. Suppose also that the homogeneous part of degree d of f has an absolutely irreducible non-linear factor in $\mathbf{Q}[x_1,..., x_n]$. Then $n(f;B)=O_{d,n,\varepsilon}(B^{n-2+\varepsilon})$ if $d \geq 4$ and $n(f;B)=O_{n,\varepsilon}(B^{n-3+2/\sqrt{3}+\varepsilon})$ if $d=3$.*

The upper bound is essentially optimal when $d \geq 4$ as $n(f;B)$ has growth order $n-2$ for polynomials of the shape $f=x_1g-x_2h$. The theorem can also be given a more geometric formulation (cf. lemma 7) in terms of the projective hypersurface defined by
$F(X_0,...,X_n):= X_0^d f(X_1/X_0,..., X_n/X_0)$.

**Theorem 1** (*alternative formulation*)  *Let $X \subset \mathbf{P}^n$ be an integral hypersurface defined by a homogeneous polynomial $F(X_0,...,X_n) \in \mathbf{Q}[X_0,..., X_n]$ of degree d and $H_0 \subset \mathbf{P}^n$ be the hyperplane defined by $X_0=0$. Suppose that $H_0$ does not contain a projective linear space of dimension $n-3$ of points of multiplicity d on X and suppose also that there is a geometrically irreducible component of degree at least two on $X \cap H_0$. Let $f(x_1,..., x_n)=F(1, x_1,..., x_n)$. Then $n(f;B)=O_{d,n,\varepsilon}(B^{n-2+\varepsilon})$ if $d \geq 4$ and $n(f;B)=O_{n,\varepsilon}(B^{n-3+2/\sqrt{3}+\varepsilon})$ if $d=3$.*

In case $X \cap H_0$ is geometrically integral, then theorem 1 was first proved by Browning, Heath-Brown and the author in [BHS] for $d \geq 6$ and after that for $d \geq 3$ in [$S_3$]. We then announced theorem 1 in a talk 2010 just after the first preprint version of [$S_3$] appeared. The proof of this theorem for surfaces is to a large extent a reexamination of the proofs in [BHS] and [$S_3$]. We then deduce the theorem for higher-dimensionsal varieties by means of repeated summation over hyperplane sections. This part of the proof is more complicated than in [BHS], but based on techniques used in previous papers of the author.

Theorem 1 was recently rediscovered by Vermeulen [V]. His result is more general than our theorem as he does not demand the existence of a component of degree at least two on $X \cap H_0$. But one cannot remove the other hypothesis as $n(f;B) \gg B^{n-2+1/d}$ for $f(x_1,...,x_n)= x_1- x_2^d$.

We shall in this note also treat affine quadrics, which were not studied in [BHS] or [$S_3$].

**Theorem 2** *Let $f(x_1,\ldots, x_n) \in \mathbf{Q}[x_1,\ldots, x_n]$ be an irreducible quadratic polynomial in $n \geq 3$ variables. Suppose that f cannot be expressed as a polynomial in two linear forms or that the homogeneous quadratic part of f is irreducible in $\mathbf{Q}[x_1,\ldots, x_n]$. Then $n(f;B) = O_{n,\varepsilon}(B^{n-2+\varepsilon})$.*

To prove theorem 2, we use again hyperplane sections to reduce to the case where $n=3$. But the proof of this case is very different from the proof of theorem 1 for surfaces where results from [$S_3$] play an essential role. We will instead use a refinement of a lemma of Heath-Brown [$H_1$] on integral points of certain conics.

We shall in the proofs of the theorems in the following sections assume that $f$ is absolutely irreducible as we have the following known lemma.

**Lemma 1** *Let $f(x_1,\ldots, x_n) \in \mathbf{Q}[x_1,\ldots, x_n]$ be an irreducible polynomial of degree d, which is not absolutely irreducible. Then $n(f;B) = O_{d,n}(B^{n-2})$.*

*Proof.* Let $X \subset \mathbf{P}^n$ be the hypersurface defined by $F(X_0,\ldots,X_n) := X_0^d f(X_1/X_0,\ldots, X_n/X_0)$. Then $X$ is integral but not geometrically integral. There exists therefore by the proof of [$S_1$, thm 2.1] a set of $O_d(1)$ proper subvarieties $Y_j$, $j \in J$ of $X$ of degree $O_{d,n}(1)$ with $X(\mathbf{Q}) = \cup_j Y_j(\mathbf{Q})$. It is thus enough to show that there are only $O_{\delta,n}(B^{n-2})$ integral $n$-tuples $(a_1,\ldots,a_n) \in [-B, B]^n$ with $(1,a_1,\ldots,a_n)$ on a subvariety $Y \subset \mathbf{P}^n$ of codimension 2 of degree $\delta$. If $X_n - bX_0$ vanishes on $Y$, then we use the trivial estimate $n(f_b;B) = O_{d,n}(B^{n-2})$ for $f_b(x_1,\ldots,x_{n-1}) = f(x_1,\ldots,x_{n-1},b)$. Otherwise, we use the induction hypothesis for all intersections of $Y$ with hyperplanes defined by $X_n - bX_0$ for integers $b$ in $[-B, B]$.

*Acknowledgement*: I would like to thank Raf Cluckers for informing me about Floris Vermeulen's work and for bringing me in contact with him.

## 2 Affine surfaces of degree at least three

We shall in this section prove theorem 1 for surfaces. We will thereby make essential use of the following deep result in section 7 of [$S_3$], which was proved by means of the authors global determinant method.

**Lemma 2** *Let $X \subset \mathbf{P}^3$ be a geometrically integral surface of degree d and $B \geq 1$. Then there exists a set S of $O_{d,\varepsilon}(B^{1/\sqrt{d}+\varepsilon})$ curves of degrees bounded solely in terms of d such that all but $O_{d,\varepsilon}(B^{2/\sqrt{d}+\varepsilon})$ points of the form $(1, a_1, a_2, a_3)$ with $(a_1, a_2, a_3)$ in $[-B, B]^3 \cap \mathbf{Z}^3$ lie on one of these curves.*

To count integral points on the curves that appear in lemma 2, we use the following result of Pila [P].

**Lemma 3** *Let $C \subset \mathbf{P}^3$ be an integral curve of of degree $\delta$ defined over $\mathbf{Q}$ and $B \geq 1$. There are then $O_{\delta,\varepsilon}(B^{1/\delta+\varepsilon})$ integral triples $(a_1, a_2, a_3)$ in $[-B, B]^3$ with $(1, a_1, a_2, a_3)$ in $C(\mathbf{Q})$.*

We shall also need the following two geometric lemmas to control the contribution from the lines in $S$. The first is similar to lemma 9 in [BHS].

**Lemma 4** *Let $X \subset \mathbf{P}^3$ be a geometrically integral surface of degree d and $P \in X$. Suppose that X is not covered by lines on X through P. There are then $O_d(1)$ lines on X passing through P.*

*Proof.* Let $H$ be the projective space parameterising (possibly reducible or non-reduced) surfaces $X \subset \mathbf{P}^3$ of degree $d$ and $W$ be the closed subset of $\mathbf{P}^3 \times \mathbf{P}^3 \times H$ representing triples $(P,Q,X)$ for which there is a line $l \subset X$ passing through $P$ and $Q$. By specialising the set of trihomogeneous polynomials defining $W$ to the given $P$ and $X$ we then obtain that the union $W_{X,P} \subset \mathbf{P}^3$ of the lines on $X$ passing through $P$ is defined by $O_d(1)$ polynomials of degree $O_d(1)$. This proves the assertion as $W_{X,P} \neq \mathbf{P}^3$ by the hypothesis.

**Lemma 5** *Let $X \subset \mathbf{P}^3$ be a geometrically integral surface of degree $d \geq 2$ with an absolutely irreducible component $C$ of degree at least two on $X \cap H_0$ for some plane $H_0 \subset \mathbf{P}^3$. There are then for each line $l_0$ on $X \cap H_0$ only $O_d(1)$ lines on $X$ intersecting $l_0$.*

*Proof.* Let $\mathbf{G}(2,4) \subset \mathbf{P}^5$ be the Grassmannian of lines $l \subset \mathbf{P}^3$ with its Plücker embedding and $H$ be the projective space parameterising surfaces $X \subset \mathbf{P}^3$ of degree $d$. Further, let $I$ be the closed subscheme of $\mathbf{G}(2,4) \times H$ parameterising pairs $(l, X)$ where $l \subset X$ and $l \cap l_0 \neq \emptyset$ and $F \subset \mathbf{P}^3 \times I$ be the associated family of lines in $\mathbf{P}^3$. We have then by the main theorem of elimination theory that the projection of $F \subset \mathbf{P}^3 \times \mathbf{G}(2,4) \times H$ on $\mathbf{P}^3 \times H$ is a closed subscheme $T$ of $\mathbf{P}^3 \times H$. The fibre $T_X \subset \mathbf{P}^3$ of the projection $\mathrm{pr}_2: T \to H$ over the point representing $X \subset \mathbf{P}^3$ will thus be defined by a set $M$ of homogeneous polynomials of degree $O_d(1)$. But the underlying set of $T_X$ is the union of the lines on $X$ intersecting $l_0$. We have therefore that $M \neq \emptyset$ as $C \not\subset T_X$, thereby finishing the proof.

We shall also make use of the following result (see (6.7) in [BHS]) to count points on lines. For $n$-tuples $\mathbf{a}=(a_1,...,a_n) \in \mathbf{Z}^n$, we will write $|\mathbf{a}|$ for $\max(|a_1|,...,|a_n|)$

**Lemma 6** *Let $l \subset \mathbf{P}^n$ be a line which intersects $H_0 \subset \mathbf{P}^n$ in a point respresented by the primitive $n$-tuple $\mathbf{b}=(b_1,...,b_n) \in \mathbf{Z}^n$. There are then $O(1+B/|\mathbf{b}|)$ $n$-tuples $\mathbf{a}=(a_1,...,a_n) \in \mathbf{Z}^n$ with $|\mathbf{a}| \leq B$.*

We are now in a position to prove theorem 1 when $n=3$.

*Proof.* It is by lemma 2 enough to consider the contribution to $n(f;B)$ from a set $S$ of $O_{d,\varepsilon}(B^{1/\sqrt{d}+\varepsilon})$ curves of degrees bounded solely in terms of $d$. Let $S_1$ be the set of curves of degree at least two in $S$. We have then by lemma 3 that each of these curves contribute with $O_{d,\varepsilon}(B^{1/2+\varepsilon})$. The total contribution to $n(f;B)$ from the curves in $S_1$ is thus $O_{d,\varepsilon}(B^{1/\sqrt{d}+1/2+\varepsilon})$, which is acceptable. Next, let $S_2$ be the set of lines in $S$, which intersect $X \cap H_0$ in components of degree $\geq 2$. We may then by lemma 4 and lemma 6 apply the argument with dyadic summation in section 6.1 in [BHS] to conclude that the total contribution from $S_2$ to $n(f;B)$ is $O_{d,\varepsilon}(B^{1+\varepsilon})$. It thus only remains to consider the contribution from the subset $S_0$ of $S$ of lines, which intersect $X \cap H_0$ in a line $l_0$. But $\#S_0 = O_d(1)$ by lemma 5 and the contribution from each line $O_d(B)$. The total contribution from $S_0$ to $n(f;B)$ is thus $O_d(B)$, which completes the proof of theorem 1 for surfaces.

**Remark 1** We have here used the hypothesis that there is an absolutely irreducible component of degree at least two on $X \cap H_0$ to control the contribution from $S_0$. There may otherwise be infinitely many lines on $X$ intersecting a line on $X \cap H_0$. But Vermeulen [V] has recently proved that the assertion in theorem 1 holds even in the case where there exists such a family of lines. It suffices that $X$ is not a cone with vertex in $H_0$.

## 3 Affine hypersurfaces of degree at least three

The aim of this section is to prove theorem 1 by means of induction with respect to $n$. We will for this need information about the geometry of the hyperplane sections of the hypersurface.

**Lemma 7** *Let $X \subset \mathbf{P}^n$ be a hypersurface of degree $d$ (i.e. a closed subscheme defined by an arbitrary homogeneous polynomial $F(X_0,...,X_n)$ of degree $d$).*
(*a*) *A point $P$ on $X$ is of multiplicity $d$ if and only if any other point of $X$ lies on a line on $X$ passing through $P$.*
(*b*) *The points of multiplicity $d$ on $X$ form a (possible empty) projective linear subspace $\Lambda(X)$.*
(*c*) *Let $H \subset \mathbf{P}^n$ be a hyperplane not containing $\Lambda(X)$. Then $\Lambda(X \cap H) = \Lambda(X) \cap H$.*

*Proof.* (*a*) It is enough to consider the Taylor expansion in the case $P=(1,0,..,0)$.

(*b*) Suppose $n \geq 2$ and let $R$ be a point on the line $L$ through two points $P$ and $Q$ of multiplicity $d$ on $X$. It is then enough to prove that any point $S \neq R$ of $X$ lies on a line on $X$ passing through $R$. We may after a coordinate change assume that $X_3 = ... = X_n = 0$ at these four points and hence also that $n=2$ with $P=(1,0,0)$ and $Q=(0,1,0)$. But then $F(X_0, X_1, X_2) = cX_2^d$ and we are done.

(*c*) It is trivial that $\Lambda(X) \cap H \subseteq \Lambda(X \cap H)$. To see that $\Lambda(X \cap H) \subseteq \Lambda(X)$, fix $R \in \Lambda(X) \setminus H$ and let $P \in \Lambda(X \cap H)$ and $Q \in X \setminus P$. We must then show that the line $l$ between $P$ and $Q$ lies on $X$. This is clear if $R \in l$. So suppose that $P$, $Q$ and $R$ span a plane $\Pi$ and let $S \in H$ be the point where the line between $R \in \Lambda(X)$ and $Q \in X$ intersects $H$. Then $S \in X$, which together with $P \in \Lambda(X \cap H)$ implies that the line $\Pi \cap H$ between these two points lies on $X$. But then $\Pi \subseteq X$ as all points on the lines between $R \in \Lambda(X) \setminus H$ and $\Pi \cap H \subseteq X$ belong to $X$. Hence $l \subseteq X$, as asserted.

The proof of the following lemma is similar to proofs in my previous papers (cf. e.g. lemma 9 in [BHS] and results in [S$_2$]).

**Lemma 8** *Let $X \subset \mathbf{P}^n$ be a geometrically integral hypersurface of degree $d \geq 2$ and $H_0 \subset \mathbf{P}^n$ be a hyperplane with $\Lambda(X) \cap H_0 \neq \emptyset$. Let $\mathbf{P}^{n\vee}$ be the dual projective space and $S_X \subset \mathbf{P}^{n\vee}$ be the subset which parameterises hyperplanes $H$ with $\dim \Lambda(X \cap H) \cap H_0 \geq \dim \Lambda(X) \cap H_0$. Then $S_X$ is a closed proper subset of $\mathbf{P}^{n\vee}$ defined by $O_{d,n}(1)$ homogeneous polynomials of degree $O_{d,n}(1)$.*

*Proof.* Let $\mathcal{H}$ be the projective space parameterising hypersurfaces $X$ of degree $d$, $\mathcal{F} \subset \mathbf{P}^n \times \mathcal{H}$ be the universal family of hypersurfaces of degree $d$, and $\mathcal{F}_0 \subset H_0 \times \mathcal{P}$ the intersection of $\mathcal{F}$ with $H_0 \times \mathcal{H}$ in $\mathbf{P}^n \times \mathcal{H}$. Further, let $I \subset \mathcal{F}_0 \times \mathbf{P}^{n\vee}$ be the closed subset (cf. remark 2) of points $(P, X, H)$ in $H_0 \times \mathcal{H} \times \mathbf{P}^{n\vee}$ such that there is an $(n-3)$-dimensional family of lines on $X \cap H$ through $P$ and $\Sigma_r$ be the set of all $(X, H) \in \mathcal{H} \times \mathbf{P}^{n\vee}$ for which the fibre of the projection $I \to \mathcal{H} \times \mathbf{P}^{n\vee}$ over $(X, H)$ is of dimension at least $r$. Then $\Sigma_r$ is a closed subset of $\mathcal{H} \times \mathbf{P}^{n\vee}$ by Chevalley's upper semi-continuity theorem (see EGA IV 13.1.5) for proper morphisms. If we fix a geometrically integral hypersurface $X$ of degree $d$ and let $r = \dim \Lambda(X) \cap H_0$, we conclude in particular that the fibre $\Sigma_X \subset \mathbf{P}^{n\vee}$ of $\mathrm{pr}_1: \Sigma_r \to \mathcal{H}$ over $X \in \mathcal{H}$ is defined by $O_{d,n}(1)$ homogeneous polynomials of degree $O_{d,n}(1)$. The points on this fibre $\Sigma_X$ will by lemma 7(*a*) represent the hyperplanes

$H \subset \mathbf{P}^n$ such that dim $\Lambda(X \cap H) \cap H_0 \geq r$. The underlying set of $\Sigma_X$ is thus $S_X$. We have further by lemma 7(c) that $H \notin S_X$ for $H$ not containing $\Lambda(X) \cap H_0$. Hence $S_X \neq \mathbf{P}^{n\vee}$, which finishes the proof.

**Remark 2** Let $\mathcal{G}$ be the Grassmannian of lines $l$ in $\mathbf{P}^n$ and $J$ the closed set of all quadruples $(l, P, X, H)$ in $\mathcal{G} \times H_0 \times \mathcal{H} \times \mathbf{P}^{n\vee}$ such that $l \subseteq X \cap H$ and $P \in l \cap H_0$. Then $I$ is the set of points on $\mathcal{F}_0 \times \mathbf{P}^{n\vee}$ where the fibre of the projection from $J$ to $\mathcal{F}_0 \times \mathbf{P}^{n\vee}$ is of dimension at least $n-3$. Hence $I$ is closed in $\mathcal{F}_0 \times \mathbf{P}^{n\vee}$ by the semi-continuity theorem quoted above.

**Lemma 9** *Let $n \geq 4$ and $X$ be a geometrically integral hypersurface of degree $d \geq 2$ in $\mathbf{P}^n$ and $H_0 \subset \mathbf{P}^n$ be a hyperplane satisifying the following conditions.*

*(i) There is a geometrically irreducible component of degree at least two on $X \cap H_0$.*
*(ii) $H_0$ does not contain a projective linear space of codimension two where all points are of multiplicity $d$ on $X$.*

*There exists then a hypersurface $Y$ of degree $O_{d,n}(1)$ in the dual projective space $\mathbf{P}^{n\vee}$ such that for every hyperplane $H \in \mathbf{P}^{n\vee} \setminus Y$ we have that $X \cap H$ is geometrically integral, $H_0 \cap H \neq H$ and such that the pair $(X^*, H_0^*) = (X \cap H, H_0 \cap H)$ satisfies (i) and (ii). We may also choose $H$ to be defined over $\mathbf{Q}$ if $X$ and $H_0$ are defined over $\mathbf{Q}$.*

*Proof.* It follows from [BrS, lemma 2.2.1] that there exists a hypersurface $W$ of degree $O_{d,n}(1)$ in $\mathbf{P}^{n\vee}$ such that $X \cap H$ is geometrically integral for $H \in \mathbf{P}^{n\vee} \setminus W$. We may also as $n \geq 4$ apply the result in (op.cit.) to the geometrically integral hypersurface in $H_0$ given by a geometrically irreducible component $V$ of degree $\geq 2$ of $X \cap H_0$. We then obtain a hypersurface $Z_0$ of degree $O_{d,n}(1)$ in the dual projective space $H_0^\vee$ of $H_0$ such that all intersections of $V \cap H'$ in $H_0$ with hyperplanes $H' \in H_0^\vee \setminus Z_0 \subset H_0$ are geometrically integral. Now let $\pi: \mathbf{P}^{n\vee} \setminus P_0 \to H_0^\vee$ be the projection from the point $P_0 \in \mathbf{P}^{n\vee}$ corresponding to $H_0 \subset \mathbf{P}^n$ and $Z$ be the closure of $\pi^{-1}(Z_0)$ in $\mathbf{P}^{n\vee}$. Then $Z$ is a hypersurface in $\mathbf{P}^{n\vee}$ of degree $O_{d,n}(1)$ such that $V \cap H$ is a geometrically irreducible component of degree at least two of $X^* \cap H_0^* = X \cap H_0 \cap H$ for all $H \in \mathbf{P}^{n\vee} \setminus Z$. Finally, if (ii) holds for $(X, H_0)$, then dim $\Lambda(X) \cap H_0 \leq$ dim $H_0 - 3$ (cf. lemma 7(b)). There exists thus by lemma 8 a hypersurface $Y$ of degree $O_{d,n}(1)$ in $\mathbf{P}^{n\vee}$ such that $\Lambda(X^* \cap H_0^*) = \Lambda(X \cap H) \cap H_0$ is of dimension at most dim $H_0 \cap H - 3$ for all $H \in \mathbf{P}^{n\vee} \setminus Y$. Hence (ii) holds for those $(X \cap H, H_0 \cap H)$ and we are done.

We are now in a position to prove theorem 1 in all dimensions by means of induction with respect to $n$.

*Proof.* We have already shown the theorem when $n=3$ in section 2. So suppose that $n \geq 4$ and let $Y \subset \mathbf{P}^{n\vee}$ be as in lemma 9. There are then only $O_{d,n}(B^n)$ linear forms $L(x_0, \ldots, x_n)$ over $\mathbf{Z}$ with $\|L\| \leq B$ which give rise to hyperplanes $H \subset \mathbf{P}^n$ in $Y$. We may therefore just as in the proof of lemma 8 in [BHS] find a linear form $L(X_0, \ldots, X_n)$ over $\mathbf{Z}$ with $\|L\|$ bounded solely in terms of $d$ and $n$ such that $L$ defines a hyperplane $H \in \mathbf{P}^{n\vee} \setminus Y$. We may also after a coordinate change assume that this holds for $L = X_n$. Let $H_b \subset \mathbf{P}^n$ be the hyperplane defined by $X_n - bX_0 = 0$ and $f_b(x_1, \ldots, x_{n-1}) = f(x_1, \ldots, x_{n-1}, b)$ for $b \in \mathbf{Z}$. Then $f_b$ will satisfy the assumptions in the theorem for $b \in \mathbf{Z}$ with $H_b \notin Y$. We have thus by the induction hypothesis that $n(f_b; B) = O_{n,\varepsilon}(B^{n-4+2/\sqrt{3}+\varepsilon})$ if $d=3$ and $n(f_b; B) = O_{d,n,\varepsilon}(B^{n-3+\varepsilon})$ if $d \geq 4$ for $b$ with $H_b \notin Y$. There are further by lemma 9 only $O_{d,n}(1)$ hyperplanes $H_b \in Y$ as $H_b = H \notin \mathbf{P}^{n\vee} \setminus Y$ for $b=0$ and we have for each $f_b$ the trivial estimate

$n(f_b;B)= O_{d,n}(B^{n-2})$. We therefore obtain the desired estimate of $n(f;B)$ from the above estimates of $n(f_b;B)$ for integers $b$ in $[-B, B]$.

## 4 Affine surfaces of degree two

We shall in this section prove theorem 2 for quadratic surfaces with methods different from the ones used in section 2. We will thereby write $||P||$ for the maximum modulus of the coefficients of a polynomial $P$.

**Lemma 10** *There exists an absolute constant C such that the following holds for any quadratic form $Q(x_1, x_2, x_3)$ over $\mathbf{Q}$ of rank at least two.*
*(a) There exists a linear form $L(x_1, x_2, x_3)$ over $\mathbf{Z}$ with $||L|| \leq C$ such that the line and conic defined by L and Q intersect transversally at two different points in $\mathbf{P}^2$.*
*(b) If Q is irreducible, then there exists a linear form $L(x_1, x_2, x_3)$ over $\mathbf{Z}$ with $||L|| \leq C$ such that $\mathbf{a}=(0, 0, 0)$ is the only triple in $\mathbf{Q}^3$ with $Q(\mathbf{a})=L(\mathbf{a})=0$.*

*Proof.* There are only $O(B^2)$ primitive triples of integers $\mathbf{a}$ in $[-B, B]^3$ with $Q(\mathbf{a})=0$. We may thus find a primitive triple $\mathbf{a}=(a_1, a_2, a_3)$ of integers with $Q(\mathbf{a}) \neq 0$ where $\max|a_i| \leq C_1$ for some absolute constant $C_1$ and three linearly independent linear forms $L_1, L_2, L_3 \in \mathbf{Z}[x_1, x_2, x_3]$ with uniformly bounded $||L_i||$ such that $L_1(\mathbf{a})=1$ and $L_2(\mathbf{a})=L_3(\mathbf{a})=0$. We may therefore assume that $Q(1, 0, 0) \neq 0$ and consider the morphism from the conic $Y \subset \mathbf{P}^2$ defined by $H$ to $\mathbf{P}^1$, which sends $(x_1, x_2, x_3)$ to $(x_2, x_3)$. But this morphism cannot ramify at more than two points. There will thus exist a pair $(c_2, c_3) \in \mathbf{Z}^2 \cap [-1,1]^2$ such that (a) holds for $L(x_1, x_2, x_3)=c_3x_2-c_2x_3$. To prove (b), we use instead theorem 3 in [H$_2$], which tells us that there are only $O_\varepsilon(B^{1+\varepsilon})$ primitive triples $\mathbf{c}=(c_1, c_2, c_3)$ in $[-B, B]^3$ with $Q(\mathbf{c})=0$. There exists therefore a pair $(c_2,c_3) \in \mathbf{Z}^2$ with $\max|c_i| \leq C$ for some absolute constant $C$ such that $(c_2, c_3) \in \mathbf{P}^1$ is not the image of a rational point on $Y$. Then $L(x_1, x_2, x_3)=c_3x_2-c_2x_3$ will be a linear form with the desired properties.

The following result is a generalisation of theorem 3 in [H$_1$].

**Lemma 11** *Let $f(x_1, x_2, x_3)$ be a quadratic polynomial with coefficients in $\mathbf{Z}$ such that the homogeneous quadratic part $Q_0(x_1, x_2)=c_{11}x_1^2+c_{12}x_1x_2+c_{22}x_2^2$ of $f(x_1, x_2, 0)$ is of rank two. Then the following holds.*
*(a) Suppose that f is not a polynomial in two linear forms. There are then for all but at most two $k \in [-B, B] \cap \mathbf{Z}$ only $O_\varepsilon((||f||B)^\varepsilon)$ integral pairs $(x_1, x_2)$ in $[-B, B]^2$ with $f(x_1, x_2, k)=0$.*
*(b) Suppose that $Q_0$ is anisotropic. There are then for all $k \in [-B, B] \cap \mathbf{Z}$ only $O_\varepsilon((||f||B)^\varepsilon)$ integral pairs in $[-B, B]^2$ with $f(x_1, x_2, k)=0$.*

*Proof.* We may assume that $c_{11} \neq 0$ after permuting $x_1, x_2$ if necessary. Let $x_1^*=2c_{11}x_1+c_{12}x_2$. We may then view $4c_{11}^2 f$ as a polynomial $g(x_1^*, x_2, x_3)$ in $\mathbf{Z}[x_1^*, x_2, x_3]$ with $||g||=O(||f||^3)$ and with $c_{11}x_1^{*2}+(4c_{11}^2c_{22}-c_{11}c_{12}^2)x_{22}^2$ as leading form of $g(x_1^*, x_2, 0)$. Since $f(x_1, x_2, k)=0$ is equivalent to $g(x_1^*, x_2, k)=0$ and $x_1^*=O(||f||B))$ for $(x_1, x_2) \in [-B, B]^2$, we have thus reduced to the case where $c_{12}=0$ and we have then that $c_{11} \neq 0$ and $c_{22} \neq 0$ by the hypothesis on $Q_0$. Now let $f(x_1, x_2, x_3)=c_{11}x_1^2+c_{13}x_1x_3+c_{22}x_2^2+c_{23}x_2x_3+c_1x_1+c_2x_2+c_{33}x_3^2+c_3x_3+c$ and $x_i^*=2c_{ii}x_i+c_{i3}x_3$ for $i=1,2$. Then $4c_{11}^2c_{22}^2 f(x_1, x_2, x_3)= c_{11}c_{22}^2(x_1^*+c_1)^2+ c_{11}^2c_{22}(x_2^*+c_2)^2+q(x_3)$ for some $q(x_3)$ in $\mathbf{Z}[x_3]$. It is also clear that the hypothesis is preserved under the transition from $f(x_1, x_2, x_3)$ to $g(x_1^*, x_2, x_3)=4c_{11}^2c_{22}^2 f \in \mathbf{Z}[x_1^*, x_2^*, x_3]$ and that we have a uniform bound $||g||=O((||f||)^{O(1)})$. We have thus reduced to the case where $f$ has the shape $a_1(x_1+b_1x_3)^2+a_2(x_2+b_2x_3)^2+q(x_3)$ for a quadratic polynomial $q(x_3) \in \mathbf{Z}[x_3]$. It is then known that there are only $O_\varepsilon((\|f\|B)^\varepsilon)$ integral pairs $(x_1, x_2)$ in $[-B, B]^2$ with $f(x_1, x_2, k)=0$ for $k \in [-B, B] \cap \mathbf{Z}$ with $q(k) \neq 0$ as explained in the proof of theorem 3 in [H$_1$]. This proves the first assertion as $q$ cannot be the zero polynomial under the hypothesis in (*a*). To deduce (*b*), note that $x_1+b_1k=x_2+b_2k=0$ if $q(k)=0$ and $Q_0$ is anisotropic.

**Lemma 12** *Let $f(x_1, x_2, x_3) \in \mathbf{Z}[x_1, x_2, x_3]$ be a polynomial of degree two with coprime coefficients and $B \geq 1$. Then one of the following holds.*

(*i*) $\|f\| = O(B^{20})$.
(*ii*) *There exists another quadratic polynomial $g(x_1, x_2, x_3) \in \mathbf{Z}[x_1, x_2, x_3]$ not proportional to $f$ such that $g(\mathbf{a})=0$ for all $\mathbf{a}=(a_1, a_2, a_3)$ in $[-B, B]^3 \cap \mathbf{Z}^3$ with $f(\mathbf{a})=0$.*

*Proof.* Let $G \in \mathbf{Z}[X_0, X_1, X_2, X_3]$ be the quadratic form with $G(1, x_1, x_2, x_3)=f$. The result is then a special case of lemma 5 in [BHS] applied to $G$.

We are now in a position to prove theorem 2 when $n=3$.

*Proof.* We apply lemma 10 to the homogeneous quadratic part $Q$ of $f$. We may then assume that $c_3 \neq 0$ in the linear form $L=c_1x_1+c_2x_2+c_3x_3$ and express $c_3^2 f(x_1, x_2, x_3)$ as a polynomial in $\mathbf{Z}[x_1, x_2, x_3^*]$ for $x_3^* = c_1x_1+c_2x_2+c_3x_3$. It is hence enough to prove the theorem in the case where the assertions of lemma 10 hold for $L=x_3$ so that $Q_0(x_1, x_2)=Q(x_1, x_2, 0)$ is of rank two. Suppose first that $f$ is not a polynomial in two linear forms. We obtain then from lemma 11(*a*) that there are $O_\varepsilon((\|f\|B)^\varepsilon)$ integral pairs $(x_1, x_2)$ in $[-B, B]^2$ with $f(x_1, x_2, k)=0$ for all but at most two $k \in [-B, B] \cap \mathbf{Z}$. We have further $O(B)$ such pairs $(x_1, x_2)$ in $[-B, B]^2$ for the exceptional values of $k$. We obtain hence after summing over $k$ that there are $O_\varepsilon(B^{1+\varepsilon}\|f\|^\varepsilon)$ integral triples $\mathbf{a} \in [-B, B]^3$ with $f(\mathbf{a})=0$. If instead $Q$ is irreducible over $\mathbf{Q}$, then $Q_0(x_1, x_2)$ is anisotropic by lemma 10(*b*). We obtain thus the same bound $O_\varepsilon(B^{1+\varepsilon}\|f\|^\varepsilon)$ as before by lemma 11(*b*). We have therefore proved the theorem if $\|f\| = O(B^{20})$ as in case (*i*) of lemma 12. If instead we are in case (*ii*), then we obtain the result from lemma 3 or by more elementary arguments, therereby completing the proof.

## 5 Affine hypersurfaces of degree two

We now prove theorem 2 in all dimensions by means of induction with respect to $n$. We shall for this use the following lemma, which follows from the proof of lemma 9.

**Lemma 13** *Let $n \geq 4$ and $X$ be a geometrically integral hypersurface of degree $d \geq 2$ in $\mathbf{P}^n$ and $H_0 \subset \mathbf{P}^n$ be a hyperplane, which does not contain a projective linear space of codimension two where all points are of multiplicity $d$ on $X$. There exists then a hypersurface $Y$ of degree $O_{d,n}(1)$ in the dual projective space $\mathbf{P}^{n\vee}$ such that for every hyperplane $H \in \mathbf{P}^{n\vee} \setminus Y$, we have that $X \cap H$ is geometrically integral, $H_0 \cap H \neq H$ and such that the pair $(X^*, H_0^*) = (X \cap H, H_0 \cap H)$ satisfies assertion (i) in lemma 9. We may also choose this hypersurface $Y$ to be defined over $\mathbf{Q}$ if $X$ and $H_0$ are defined over $\mathbf{Q}$.*

*Proof* (of theorem 2). We have already shown the theorem when $n=3$ in section 4. So suppose that $n \geq 4$ and let us first consider the case where $f$ cannot be expressed as a polynomial in two linear forms, which is equivalent to that $H_0$ does not contain a projective linear space of dimenion $n-3$ of points of multiplicity $d$ on $X$. The proof of theorem 1 in section 3 will then extend almost verbatim to this case if we just use lemma 13 instead of lemma 12. It thus only remains to treat the case where the quadratic part $Q$ of $f$ is irreducible over **Q**. We may also assume that rank $Q=2$ as $f$ cannot be a polynomial in two linear forms if rank $Q \geq 3$. Now let $Z \subset \mathbf{P}^{n-1}$ be the quadric defined by $Q$ and $\Lambda(Z)$ be its singular locus. $\Lambda(Z)$ is then a projective linear subspace of codimension two in $\mathbf{P}^{n-1}$ corresponding to a line in the dual projective space $\mathbf{P}^{(n-1)\vee}$. There are thus $O_n(B^2)$ linear forms $L(x_0,\ldots, x_n)$ over **Z** with $\|L\| \leq B$ and $L=0$ on $\Lambda(Z)$. We may therefore find a linear form $L(x_0,\ldots, x_n)$ over **Z** with $\|L\|$ bounded solely in terms of $n$ such that $L \neq 0$ on $\Lambda(Z)$ and we may further after a linear coordinate change assume that $L=x_n$. The quadratic form $q(x_1,\ldots, x_{n-1}) = Q(x_1,\ldots, x_{n-1},0)$ is then of rank two and irreducible over **Q**. As $q$ is also the homogeneous quadratic part of all $f_b=f(x_1,\ldots, x_{n-1},b)$, we have thus by the induction hypothesis that $n(f_b;B) = O_{n,\varepsilon}(B^{n-3+\varepsilon})$ for all $b \in \mathbf{Z}$. But then $n(f;B)=O_{n,\varepsilon}(B^{n-2+\varepsilon})$ as $n(f;B)$ is the sum of all $n(f_b;B)$ for $b \in [-B, B] \cap \mathbf{Z}$.